\documentclass[12pt]{amsart}
\usepackage{amsmath}
\usepackage{amssymb}
\usepackage[latin1]{inputenc}
\usepackage[T1]{fontenc}
\newtheorem{lemma}{Lemma}

\newtheorem{remark}{Remark}

\newtheorem{theorem}{Theorem}

\newtheorem{proposition}{Proposition}

\newcommand{\bra}{\langle}
\newcommand{\ket}{\rangle}
\newcommand{\vp}{\varphi}

\newcommand{\C}{\mathbb{C}}
\newcommand{\F}{\mathbb{F}}

\newcommand{\N}{\mathbb{N}}

\newcommand{\R}{\mathbb{R}}
\newcommand{\Z}{\mathbb{Z}}
\newcommand{\be}{\begin{equation}}
\newcommand{\eeq}{\end{equation}}
\newcommand{\bet}{\begin{equation*}}
\newcommand{\eeqt}{\end{equation*}}
\newcommand{\bea}{\begin{eqnarray}}
\newcommand{\eeqa}{\end{eqnarray}}
\newcommand{\beat}{\begin{eqnarray*}}
\newcommand{\eeqat}{\end{eqnarray*}}

\newcommand{\h}[1]{\mathcal{#1}}
\newcommand{\hil}{\mathcal{H}}

\newcommand{\hM}{\mathcal{M}}

\newcommand{\hT}{\mathcal{T}}

\newcommand{\V}{\mathrm{V}}

\newcommand{\tr}{{\rm tr}}
\newcommand{\lin}{{\rm lin}}
\newcommand{\rank}{{\rm rank}\,}

\newcommand{\cal}{\mathcal}

\newcommand{\hi}{\mathcal H}
\newcommand{\kb}[2]{\left|#1\right\rangle\left\langle#2\right|}

\def\<{\langle}
\def\>{\rangle}

\begin{document}
\title{A note on infinite extreme correlation matrices}

\author{J.\ Kiukas$^1$}
\author{J.-P.\ Pellonpää$^2$}

\maketitle

{\sc\small ${}^{1,2}$Department of Physics, University of Turku, FIN-20014 Turku, Finland
\it\small

$^1$e-mail: jukka.kiukas@utu.fi

$^2$e-mail: juhpello@utu.fi}

\begin{abstract}
We give a characterization for the extreme points of the convex set of correlation matrices with a countable index set. A Hermitian matrix is called a correlation matrix if it is positive semidefinite with unit diagonal entries. Using the characterization we show that there exist extreme points of any rank.
\end{abstract}
\vspace{5mm}
\newline


\section{Introduction}
Let $\cal I$ be a fixed index set, and let $\F$ denote the field of either real or complex numbers.
A mapping $C:\,{\cal I}\times{\cal I}\to\F$ is an {\em (auto)correlation function} if it is 
positive semidefinite, Hermitian (or symmetric when $\F=\R$), and $C(i,i)=1$ for all $i\in\cal I$. 
This term comes from the theory of stochastic processes.
It the case where $\cal I$ is either finite or countable, such a function is called a {\em correlation matrix}.
The set of correlation functions is convex. The problem of determining the extreme points of this set has been studied extensively, but the research is mostly concentrated on finite matrices \cite{ChVe,Lo,GrPiWa,LiTa}.

Christensen and Vesterstr\o m \cite{ChVe} considered the complex case and proved that, 
when the cardinality of $\cal I$ is greater than 3, there exists a rank 2 extreme point.
Moreover, they showed that when $\cal I$ is finite with $n$ elements,
the rank of any extreme point is at most $\sqrt{n}$. Loewy \cite{Lo} proved that, in fact, for any $r\le\sqrt{n}$ there exists a rank $r$ extreme point. Grone, Pierce, and Watkins \cite{GrPiWa} showed that, in the real case, a similar result holds but the criterion for the rank is then $r^2+r\le 2n$. 
Finally, Li and Tam were able to give a simple characterization for the extreme points in the case where $\cal I$ is finite \cite[Theorem 1 (b)]{LiTa}. This result holds in both real and complex cases.

The complex correlation matrices are used also in quantum mechanics where they appear, e.g., as the structure matrices of certain $\mathbb T$-covariant observables \cite{Ho,Pellonpää,DA}:
Finite correlation matrices occur in the context of angle observables while
the case of countable $\cal I$ is associated with
covariant phase observables (when $\cal I=\{0,1,2,...\}$) and with box localization observables (when $\cal I=\Z$). In these cases, the extreme correlation matrices correspond to the extreme $\mathbb T$-covariant measurements; accordingly, Holevo stressed the importance of characterizing them \cite{Ho}. 
For finite $\cal I$, the solution is provided by Li and Tam as mentioned above. Their result has been generalized to other finite-dimensional covariance systems by D'Ariano \cite{DA}.

The purpose of this note is to give a characterization of extreme correlation matrices in the case of countable index set (which can be chosen to be $\N$ without restricting generality). 
We consider both real and complex cases.
The result is given by Theorem 1 below. It is a generalization of Theorem 1 (b) of \cite{LiTa}. We also show that there exist exteme matrices of any rank ($\in\N\cup\{\infty\}$).

\section{Notations and basic definitions}

Let $\h M$ be the set of infinite Hermitian (i.e.\ conjugate symmetric) $\F$-valued matrices indexed by $\N=\{1,2,3,...\}$, considered as a real linear space with respect to the usual entrywise
sum and scalar multiplication. 
Let $\V$ be the vector space over $\F$ consisting of $\F$-valued sequences $c=(c_n)_{n\in\N}$ such that $c_n\ne 0$ only for finitely many $n\in\N$. For any $k\in\N$, let $e_k\in\V$ be the sequence $(\delta_{kn})_{n\in\N}$ where $\delta_{kn}$ is the Kronecker delta (so that
$(e_k)_{k\in\N}$ forms an algebraic basis of $\V$). 
Each matrix $M=(M_{nm})_{n,m\in\N}\in \h M$ defines a Hermitian sesquilinear (or symmetric bilinear when $\F=\R$) form $M:\,\V\times\V\to \F$ via 
\bet
M(c,d) = \sum_{n,m\in\N} \overline{c_n} M_{nm} d_m, \ \ \ c,d\in \V.
\eeqt
(Notice that the above sums are finite.) Let $\h C \subset \h M$ be the convex set of Hermitian positive semidefinite matrices with
unit diagonal elements, i.e.\ the set of correlation matrices.
Recall that $C\in \h M$ is positive semidefinite if and only if $C(c,c)\geq 0$ for all $c\in \V$.

Let $\ell^1(\N)$ be the Banach space (over $\F$) consisting of $\F$-valued sequences $c=(c_n)_{n\in\N}$ such that $\sum_{n\in\N}|c_n|<\infty$. We denote the $\ell^1$-norm $\sum_{n\in\N}|c_n|$ of $c\in\ell^1(\N)$ by $\|c\|_1$. 
Let $\ell^1(\N)^\times$ denote the topological antidual of $\ell^1(\N)$, i.e.\ the linear space (over $\F$) of continuous antilinear functionals $\ell^1(\N)\to\F$.
We equip $\ell^1(\N)^\times$ with the standard operator norm
$$f\mapsto\sup\{|f(c)| \mid c\in\ell^1(\N),\,\|c\|_1\le1\}.$$
Naturally, $\ell^1(\N)^\times$ is isomorphic with $\ell^\infty(\N)=\{ d:\N\to \F\mid \sup_{n\in\N} |d_n|<\infty\}$ 
but we do not actually need this fact. 

Let $C\in \h C $.
Since $C_{nn}=1$, $n\in\N$, and the principal $2\times2$-minors of $C$ are nonnegative,
we have $\sup_{n,m\in \N}|C_{nm}|\leq 1$,
implying that
\bet
\sup_{n\in \N} \Big|\sum_{m\in\N} C_{nm} c_m\Big| \leq \|c\|_1, \ \ \  c\in \ell^1(\N).
\eeqt
Hence, $C$ defines a continuous linear map $\tilde{C}: \ell^1(\N)\to \ell^1(\N)^\times$ by
\bet
\big[\tilde{C}(c)\big](d)=\sum_{n,m\in\N} \overline{d_n} C_{nm} c_m, \ \ \ c,d\in \ell^1(\N).
\eeqt
where the double series converges absolutely. 

The \emph{rank} of $C\in\h C$ is defined by
\bet
\rank C = \dim \tilde{C}(\ell^1(\N))\in\N\cup\{\infty\},
\eeqt
i.e.\ $\rank C$ is the dimension of the linear space $\tilde{C}(\ell^1(\N))\subset \ell^1(\N)^\times$.

\section{Technical lemmas}

Before we can prove the main result of this note (Theorem \ref{maintheorem}), we need some information on the structure of a matrix $C\in \h C $. 

The matrices $C\in \h C $ have the following characterization.
Let $C\in\h M$. Then $C\in \h C $ if and only if there exists a separable Hilbert space $\hil$ over $\F$ and
a sequence of unit vectors $(\eta_n)_{n\in\N}\subset\hil$ such that $C_{nm}=\bra\eta_n|\eta_m\ket$ for all $n,m\in \N$ \cite{Pellonpää} (see also \cite[Exercise 8.7]{Pa}); here 
$\bra\,\cdot\, |\,\cdot\,\ket$ denotes the inner product of $\hil$ (linear in the second argument) and later we use the symbol $\|\,\cdot\,\|$ for the norm of $\hil$.
In the above references, this result is given only in the complex case, but the real case follows easily by 
noting that a complex Hilbert space $(\hi,\bra\,\cdot\, |\,\cdot\,\ket)$ can be interpreted as the real Hilbert space equipped with the inner product ${\rm Re}\,\bra\,\cdot\, |\,\cdot\,\ket$.
Note that $\hil$ and the sequence $(\eta_n)$ are not uniquely determined by $C$.

Let $C\in\h C $, and choose a sequence $(\eta_n)$ of unit vectors in a Hilbert space $\hil$ (over $\F$) such that
$C_{nm}=\bra\eta_n|\eta_m\ket$ for all $n,m\in \N$. These choices will remain fixed in this section.
For each $c\in \ell^1(\N)$, the series $\sum_{n\in\N}c_n\eta_n$ converges (absolutely) in $\hil$, so we can define
a linear map $\Phi:\ell^1(\N)\to\hi$ by
\be\label{fii}
\Phi(c)= \sum_{n\in\N} c_n\eta_n.
\eeq
Since $\|\Phi(c)\| \leq \|c\|_1$ the map $\Phi$ is continuous.

Next define another linear map $\Phi^*:\hil\to \ell^1(\N)^\times$ by
\bet
\big[\Phi^*(\vp)\big](c) = \sum_{n\in\N} \overline{c_n}\bra \eta_n|\vp\ket, \ \ \ \vp\in \hil, \ c\in\ell^1(\N),
\eeqt
where the series converges absolutely. The map $\Phi^*$ is continuous since
\bet
\sup\big\{|[\Phi^*(\vp)](c)|\,\big|\,c\in\ell^1(\N),\,\|c\|_1\le1\big\}\leq \|\vp\|.
\eeqt
It is clear from the above formulas that we have
\be\label{transpose}
\big[\Phi^*(\vp)\big](c) = \bra \Phi(c)|\vp\ket, \ \ \ \vp\in \hil, \ c\in \ell^1(\N).
\eeq

\begin{lemma}\label{lemma}
\begin{itemize}
\item[(a)] $\overline{\lin_\F  \{ \eta_n\mid n\in \N\}}= \overline{\Phi(\ell^1(\N))}= (\ker \Phi^*)^\perp$.
\item[(b)] $\tilde{C}=\Phi^*\Phi$.
\item[(c)] $\overline{\tilde{C}(\ell^1(\N))}= \overline{\Phi^*(\hil)}$
\item[(d)] $\rank C= \dim \Phi^*(\hil)= \dim \Phi(\ell^1(\N)) = \dim \lin_\F \{ \eta_n\mid n\in \N\}$.
\end{itemize}
\end{lemma}
\begin{proof}
Since $\eta_k = \Phi(e_k)$ by \eqref{fii}, we have $\overline{ \lin_\F  \{ \eta_n\mid n\in \N\}}\subseteq \overline{\Phi(\ell^1(\N))}$. On the other hand,
the series in \eqref{fii} converges in norm, so $\Phi(\ell^1(\N))\subseteq \overline{ \lin_\F  \{ \eta_n\mid n\in \N\}}$. This proves the first equality in (a).
It follows from \eqref{transpose} that $\Phi(\ell^1(\N)) \subseteq (\ker \Phi^*)^\perp$ and
hence $\overline{\Phi(\ell^1(\N))} \subseteq (\ker \Phi^*)^\perp$ (the orthogonal complement is closed). Assume now that
$\vp_0\in (\ker \Phi^*)^\perp\cap \Phi(\ell^1(\N))^\perp$. Then $\left[\Phi^*(\vp_0)\right](c) = \bra \Phi(c)|\vp_0\ket=0$ for all $c\in \ell^1(\N)$ by \eqref{transpose},
so $\vp_0\in \ker \Phi^*$ and hence $\vp_0=0$. Now (a) is proved.

By \eqref{fii}, the continuity of $\Phi^*$, and \eqref{transpose}, we get
\beat
\big[\Phi^*\big(\Phi(c)\big)\big](d) &=& \sum_{n\in\N} c_n\big[\Phi^*(\eta_n)\big](d) = \sum_{n\in\N} c_n\bra \Phi(d)|\eta_n\ket\\
&=& \sum_{n\in\N} \sum_{m\in\N} \overline{d_m} \bra\eta_m|\eta_n\ket c_n = \big[\tilde{C}(c)\big](d), \ \ \ c,d\in \ell^1(\N).
\eeqat
This proves (b).

To prove (c), we first note that $\tilde{C}(\ell^1(\N))\subseteq \Phi^*(\hil)$ by (b), so $\overline{\tilde{C}(\ell^1(\N))}\subseteq \overline{\Phi^*(\hil)}$.
On the other hand, by using the relation $\overline{\Phi(\ell^1(\N))} = (\ker\Phi^*)^\perp$ (see (a)) and the fact that
$\Phi^*$ is continuous, one gets
\bet
\Phi^*(\hil)= \Phi^*(\overline{\Phi(\ell^1(\N))})\subseteq \overline{\Phi^*(\Phi(\ell^1(\N)))} = \overline{\tilde{C}(\ell^1(\N))},
\eeqt
and hence $\overline{\Phi^*(\hil)}\subseteq \overline{\tilde{C}(\ell^1(\N))}$, and (c) is proved.

It remains to prove (d). By (a), we have $\overline{\lin_\F  \{ \eta_n\mid n\in \N\}}=\overline{\Phi(\ell^1(\N))}= \h K$, where $\h K=(\ker \Phi^*)^{\perp}$.
Since any finite dimensional space is closed,
$\lin_\F  \{ \eta_n\mid n\in \N\}=\dim \Phi(\ell^1(\N))=\dim \h K$. Now
$\Phi^*|_{\h K}(\h K) = \Phi^*(\hil)$ since $\Phi^*(\h K^{\perp})={0}$. In addition, the map $\Phi^*|_{\h K}:\h K\to \Phi^*(\hil)$ is a linear bijection, so
$\dim \h K = \dim \Phi^*(\hil)$. 
It follows from (c) that $\dim \tilde{C}(\ell^1(\N))= \dim \Phi^*(\hil)$,
i.e.\ $\rank C=\dim \Phi^*(\hil)$. This completes the proof of (d).
\end{proof}

\begin{remark}\rm 
From Lemma 1, one sees that, without restricting generality, we may assume that $\dim\,\hil=\rank C$ and the sequence $(\eta_n)_{n\in\N}$ spans $\hil$.
In particular, if $\rank C=r<\infty$ then we may assume that $\hil=\F^r$.
\end{remark}

For each $\vp\in \hil$, we let $|\vp\ket\bra\vp|$ denote the bounded operator $\psi\mapsto \bra\vp|\psi\ket\vp$ on $\hil$. We use the symbols 
$\h L_s(\h H)$ and $\h T_s(\h H)$ for the real Banach spaces of
selfadjoint bounded and selfadjoint trace class operators on $\hil$, respectively.
The following lemma is well known in the case of a complex Hilbert space. However, real Hilbert spaces are not so frequently used in the literature, so we give the proof here for the real case.
\begin{lemma}
The topological dual of $\h T_s(\hil)$ is isomorphic to $\h L_s(\hil)$.
\end{lemma}
\begin{proof}Assume that $\hil$ is real (i.e.\ $\F=\R$).
For each $R\in\h L_s(\hil)$, define a linear mapping $F_R:\,\h T_s(\hil)\to\R$ by $F_R(T)=\tr[RT]$.
The functional $F_R$ is continuous, since
$|F_R(T)|\le\|R\|\|T\|_{\tr}$ where $\|\,\cdot\,\|$ and $\|\,\cdot\,\|_{\tr}$ are the operator and trace class norms, respectively \cite[Theorem 4.1.4 (3)]{Li}.
Let then $F:\,\h T_s(\hil)\to\R$ be a continuous functional. Define an  $\R$-linear map $\tilde F$ on the whole trace class $\h T(\hil)$ of $\hi$ by $\tilde F(T)=\frac12 F(T+T^*)$, $T\in \h T(\hil)$.
Since $\|T^*\|_{\tr}=\|T\|_{\tr}$, it follows that
$|\tilde F(T)\|\le\|F\|\|T\|_{\tr}$ and $\tilde F$ is continuous. Hence, there exists a bounded operator $S$ such that $\tilde F(T)=\tr[ST]$, $T\in \h T(\hil)$, and $\|\tilde F\|=\|S\|$ \cite[Theorem 4.1.4 (3)]{Li}.
For any $T\in\h T_s(\hil)$ we have $F(T)=\tilde F(T)=\frac12 \tr[(S+S^*)T]$. (Note that $\tr[A^*]=\tr[A]$, $A\in\h T(\hil)$.) Define $R=\frac12(S+S^*)\in\h L_s(\hil)$ so that $F=F_R$.
The mapping $R\mapsto F_R$ is isometry. Indeed,  $\|F_R\|=\|R\|$ since $\tr\big[R\frac12(T+T^*)\big]=\tr(RT)$,  $T\in\h T(\hil)$, and $\frac12\|T+T^*\|_\tr\le1$ when $\|T\|_\tr\le1$.
\end{proof}

\section{Extreme points of $\h C $}
The following theorem characterizes the extreme points of the convex set $\h C $ of $\F$-valued correlation matrices.
The characterization is an infinite-dimensional analogue of  \cite[Theorem 1 (b)]{LiTa}. 

\begin{theorem}\label{maintheorem}
Let $C\in \h C $, and let $(\eta_n)$ be a sequence of unit vectors in a Hilbert space $\hil$ (over $\F$)
such that $C_{nm}=\bra\eta_n|\eta_m\ket$ for all $n,m\in \N$. Denote $\h K=\overline{\lin_\F  \{ \eta_n\mid n\in \N\}}$ and
\bet
\h S = \lin_{\R} \{|\eta_n\ket\bra\eta_n|\mid n\in \N\} \subseteq \h T_s(\h K).
\eeqt
Then $C$ is an extreme point of $\h C $ if and only if $\h S$ is trace norm dense in $\h T_s(\h K)$.
\end{theorem}
\begin{proof}
Define $\Phi$ and $\Phi^*$ using $(\eta_n)$ as in the preceding section so that Lemma \ref{lemma} holds for them.
Notice that $\h K= (\ker \Phi^*)^\perp=\overline{\Phi(\ell^1(\N))}\subseteq \hil$ by Lemma \ref{lemma} (a).

Assume that $\h S$ is not dense in $\h T_s(\h K)$. Since the dual of $\h T_s(\h K)$ is $\h L_s(\h K)$ by Lemma 2, it follows from the
Hahn-Banach theorem that there is  an $R\in \h L_s(\h K)$ such that $R\neq 0$ and $\tr [R T]=0$ for all $T\in \h S$. Define
$\tilde{B}:\ell^1(\N)\to\ell^1(\N)^\times$ by $\tilde{B}=\Phi^*|_{\h K}R\Phi$. Now $\tilde B\neq 0$, because $\Phi^*|_{\h K}$ is injection and
$\Phi(\ell^1(\N))$ is dense in $\h K$ by Lemma \ref{lemma} (a).
Let $B_{nm}= \big[\tilde{B}(e_m)\big](e_n)$ for all $n,m\in \N$. Then by \eqref{transpose} and \eqref{fii}, we get
\bet
B_{nm}= \big[\Phi^*(R\Phi(e_m))\big](e_n) = \bra \Phi(e_n)|R\Phi(e_m)\ket = \bra \eta_n|R\eta_m\ket, \ \ \ n,m\in \N.
\eeqt
Since $R$ is selfadjoint, it follows that the matrix $(B_{nm})$ is Hermitian. In addition, $B_{nn} = \tr[R|\eta_n\ket\bra\eta_n|]=0$ for all
$n\in \N$, since $|\eta_n\ket\bra\eta_n|\in \h S$ for all $n\in\N$. Let $\epsilon = \|R\|^{-1}$ (where $\|R\|$ is the operator norm of $R$). 
Since $R$ is a selfadjoint bounded operator,
$I\pm \epsilon R\geq 0$ where $I$ is the identity operator of $\h K$.
Hence, for all $c\in \ell^1(\N)$, one gets from \eqref{transpose} that
\bet
(C\pm \epsilon B)(c,c) = \big[(\Phi^*\Phi\pm \epsilon \Phi^*R\Phi)(c)\big](c) 
= \bra\Phi(c)| (I\pm \epsilon R)\Phi(c)\ket\geq 0.
\eeqt
It follows that the matrix $C\pm \epsilon B$ is positive semidefinite. Since, in addition, $B_{nn}=0$ for all $n$, we get
$C\pm \epsilon B\in \h C $. Clearly $C= \frac 12 (C-\epsilon B)+\frac 12 (C+\epsilon B)$, so $C$ is not an extreme point.

Assume then that $C$ is not an extreme point. Then $C=\frac 12(C_1+C_2)$ for some $C_1,C_2\in \h C $, where $C_1\neq C_2$.
Define $B = \frac 12(C_1-C_2)$. Then $B\neq 0$, and $B\in \h M$, with $B_{nn} = 0$ for all $n\in \N$.
Since $C\pm B\in \h C $, one gets $C(c,c)\pm B(c,c)\geq 0$ for all $c\in \V$, and hence
\be\label{Bineq}
|B(c,c)|\leq C(c,c), \ \ \ c\in \V.
\eeq 
Next we show that
\be\label{implication}
\text{if $c\in \ker \Phi\cap \V$ and $d\in \V$ then } B(c,d) = 0.
\eeq
To that end, let $c\in \V$ be such that
$\Phi(c) = 0$, and let $d\in \V$. Then $C(c, c) =\big[\Phi^*\Phi(c)\big](c)=0$, so
also $(C+B)(c,c) = 0$ by \eqref{Bineq}. Since $C$ and $C+B$ are positive semidefinite and Hermitian, the Cauchy-Schwarz inequality gives $C(c,d)=0$ and
$(C+B)(c,d) = 0$. Hence also $B(c,d)= 0$, proving \eqref{implication}.

It follows from \eqref{implication} that $B(c_1,d_1) = B(c_2, d_2)$ whenever $c_1,c_2,d_1,d_2\in \V$ are such that
$c_1-c_1, d_1-d_2\in \ker \Phi$. (Notice that $B$ is Hermitian.) Hence, we have a well-defined 
sesquilinear form $R:\,\Phi(\V)\times\Phi(\V)\to \F$ defined by
$R(\vp,\psi) = B(c, d)$ where $\Phi(c)=\vp$ and $\Phi(d)=\psi$. 
Since $B$ is a Hermitian, also $R$ is such.

Let now $\vp\in \Phi(\V)$, and let $c\in \V$ be such that $\vp=\Phi(c)$. Then by \eqref{Bineq}, Lemma \ref{lemma} (b) and
\eqref{transpose}, we get
\bet
|R(\vp,\vp)| = |B(c,c)|\leq C(c,c) = \big[\Phi^*(\Phi(c))\big](c) = \bra \Phi(c)|\Phi(c)\ket = \|\vp\|^2,
\eeqt
so the polarization identity gives
\bet
\sup \{ |R(\vp,\psi)|\mid \vp,\psi\in \Phi(\V), \|\vp\|\leq 1,\|\psi\|\leq 1\}<\infty.
\eeqt
Since $\Phi(\V) = \lin_\F  \{ \eta_n \mid n\in \N\}$ is dense in $\h K$ by Lemma \ref{lemma} (a) it follows that there is a bounded operator
$\tilde{R}:\h K\to \h K$, such that $R(\vp,\psi) = \bra \vp|\tilde{R}\psi\ket$ for all $\vp,\psi\in \Phi(\V)$. Since $R$ is Hermitian,
$\tilde{R}\in \h L_s(\h K)$. Now $\Phi(e_n) = \eta_n$ for all $n\in \N$, so
\bet
\tr[\tilde{R}|\eta_n\ket\bra\eta_n|] = \bra \eta_n|\tilde{R}\eta_n\ket = B(e_n,e_n) = B_{nn} = 0, \ \ \ n\in \N.
\eeqt
This implies that $\tr[\tilde{R}T]=0$ for all $T\in \h S$.

Now if $\h S$ is dense in $\h T_s(\h K)$, it follows that $\tr[\tilde{R}T]=0$ for all $T\in \h T_s(\hM)$, implying $\tilde{R}=0$. But then
\bet
B_{nm} = B(e_n,e_m)=R(\eta_n,\eta_m)= \bra \eta_n |\tilde{R}\eta_m\ket = 0, \ \ \ n,m\in \N,
\eeqt
which is impossible, since $B\neq 0$. Hence, $\h S$ is not dense in $\hT_s(\h K)$. The proof is complete.
\end{proof}

The following proposition shows that there exist extreme points of $\h C $ of any rank $(\in\N\cup\{\infty\})$.
\begin{proposition}
Let $\h N$ be a (separable) Hilbert space over $\F$, and let $\{ \vp_n\mid n\in N\}\subseteq \h N\setminus\{0\}$ be a countable set dense in $\h N$. Define
\bet
C_{nm} = \|\vp_n\|^{-1}\|\vp_m\|^{-1}\bra \vp_n |\vp_m\ket, \ \ \ n,m\in \N.
\eeqt
Then the matrix $C=(C_{nm})$ is an extreme point of $\h C $ with $\rank C = \dim \h N$.
\end{proposition}
\begin{proof}
Let $\eta_n = \|\vp_n\|^{-1} \vp_n$ for all $n\in \N$. Since the set $\{ \vp_n\mid n\in \N\}$ is dense in $\h N$, it is clear that
$\overline{\lin_\F  \{ \eta_n\mid n\in \N\}} = \h N$. Therefore, $\rank C= \dim \h N$ by Lemma \ref{lemma} (d).
We proceed to show that the set $\h S = \lin_{\R} \{|\eta_n\ket\bra\eta_n|\mid n\in \N\}$ is trace norm dense in $\h T_s(\h N)$.
To that end, let $T\in \h T_s(\h N)$, $T\ne 0$, and $\epsilon >0$. Using the spectral
representation \cite{Li}, we get
\bet
T=\sum_{n\in \N} t_n |\phi_n\ket\bra\phi_n|,
\eeqt
where $t_n\in \R$ for all $n\in\N$, and $(\phi_n)_{n\in\N}$ is an orthonormal sequence in $\h N$, the series converging in the trace norm $\|\,\cdot\,\|_{\tr}$.
Choose $n_0\in \N$ such that
\be\label{eq1}
\Big\| T-\sum_{n=1}^{n_0} t_n |\phi_n\ket\bra\phi_n|\Big\|_{\tr}<\frac {\epsilon}{2}.
\eeq
Then, for each $n=0,..., n_0$, pick $k_n\in \N$ such that $\|\phi_n-\vp_{k_n}\|<\epsilon (6\|T\|_{\tr})^{-1}$. This is possible because
$\{\vp_n\mid n\in \N\}$ is dense in $\h N$. Now we use the fact that for $\psi_1,\psi_2\in \h N$, such that $\|\psi_1\|=1$ and
$\|\psi_1-\psi_2\|\leq 1$, we have
\bet
\| |\psi_1\ket\bra\psi_1|-|\psi_2\ket\bra \psi_2|\|_{\tr} \leq 3\|\psi_1-\psi_2\|.
\eeqt
This can be proved \cite[the proof of Lemma 5]{NCQM} e.g.\ by using the duality $\h C(\h N)^*=\hT(\h N)$ \cite[p.\ 60]{Li} where $\h C(\h N)$ denotes the set of compact operators.
Applying this result, we get
\be\label{eq2}
\Big\|\sum_{n=1}^{n_0} t_n |\phi_n\ket\bra\phi_n|-\sum_{n=1}^{n_0}t_n|\vp_{k_n}\ket\bra\vp_{k_n}|\Big\|_{\tr}
\leq 3\sum_{n=1}^{n_0} |t_n| \|\phi_n-\vp_{k_n}\| < \frac {\epsilon}{2}.
\eeq
Now \eqref{eq1} and \eqref{eq2} imply that
\be\label{eq3}
\Big\| T-\sum_{n=1}^{n_0} t_n |\vp_{k_n}\ket\bra\vp_{k_n}|\Big\|_{\tr}<\epsilon.
\eeq
Since $\vp_{k_n} = \|\vp_{k_n}\|\eta_{k_n}$ we have $\sum_{n=1}^{n_0} t_n |\vp_{k_n}\ket\bra\vp_{k_n}|\in \h S$. Thus $\h S$ is dense in
$\h T_s(\h N)$. Now $C$ is an extreme point by Theorem 1. The proof is complete.
\end{proof}

The following remark demonstrates that there exist extreme correlation matrices which cannot be constructed by using Proposition 1.

\begin{remark}\rm
It is easy to generalize the construction of finite rank extreme points of Li and Tam \cite[Section 2.1]{LiTa} to the infinite dimensional case. 

First we consider the complex case.
Let $r\in\N$ and choose $\hil=\C^r$. Let $(f_n)_{n=1}^r$ be an orthonormal basis of $\C^r$.
Then the basis of the finite dimensional space $\h T_s(\C^r)=\h L_s(\C^r)$ consists of operators $\kb{f_n}{f_n}$, $\kb {f_n} {f_m}+\kb {f_m} {f_n}$ and  $i\kb {f_n} {f_m}-i\kb {f_m} {f_n}$ where $n,m\in\{1,2,...,r\}$ and $n<m$.
Define $(\eta_n)$-sequence as follows. Enumerate the unit vectors $f_n$, $2^{-1/2}(f_n+f_m)$, $2^{-1/2}(f_n+if_m)$, $n,m\in\{1,...,r\},$ $n<m$, to get the first $r^2$ vectors $\eta_n$ and define $\eta_n=f_1$ for all $n>r^2$. 
Thus we have a sequence $(\eta_n)_{n\in\N}$ of unit vectors and we can define a $C\in\h C$ by $C_{nm}=\<\eta_n|\eta_m\>$. 
Since, for example, $$i\kb {f_n} {f_m}-i\kb {f_m} {f_n}=\kb{f_n}{f_n}+\kb{f_m}{f_m}-2\cdot\kb{2^{-1/2}(f_n+if_m)}{2^{-1/2}(f_n+if_m)}$$ it is easy to see that $\lin_{\R} \{|\eta_n\ket\bra\eta_n|\mid n\in \N\}=\h T_s(\C^r)$. Hence, $C$ is extreme.
Moreover, $\rank C=\dim\lin_\C\{\eta_n\mid n\in\N\}=r$. 

Let $(\tilde\eta_n)$ be any sequence of unit vectors (in some Hilbert space) such that $C_{nm}=\bra\tilde\eta_n|\tilde\eta_m\ket$. Since $\bra\tilde\eta_n|\tilde\eta_m\ket=1$ for all $n,m>r^2$ it follows that $\tilde\eta_n=\tilde\eta_{r^2+1}$ for all $n>r^2$ by Cauchy-Schwarz. Hence, there exist only finitely many different vectors $\tilde\eta_n$, so the sequence $(\tilde\eta_n)$ cannot be dense in the surface of the unit ball. It follows that the sequence $(\tilde\eta_n)$ cannot be obtained by the method of Proposition 1. 

In the real case we take $\hil=\R^r$ with an orthonormal basis $(f_n)_{n=1}^r$.
Now the operators $\kb{f_n}{f_n}$ and $\kb {f_n} {f_m}+\kb {f_m} {f_n}$, $n,m\in\{1,2,...,r\}$, $n<m$, generate $\h T_s(\R^r)$ and we can choose the first $r(r+1)/2$ $\eta_n$-vectors to be $f_n$ and $2^{-1/2}(f_n+f_m)$, $n<m$. Then we set $\eta_n=f_1$ for all $n>r(r+1)/2$. 
As before, we notice that the matrix $C_{nm}=\<\eta_n|\eta_m\>$ is an extreme point of rank $r$ and it cannot be obtained by Proposition 1.

Hence, to construct a finite rank extreme, Proposition 1 gives unnecessarily many different vectors $\eta_n$.
\end{remark}

\noindent{\bf Acknowledgments.} The authors thank Drs.\ Pekka Lahti and Kari Ylinen for fruitful discussions. One of us (J.K.) was supported by Finnish Cultural Foundation.


\begin{thebibliography}{99}

\bibitem{Pellonpää} G.\ Cassinelli, E.\ De Vito, P.\ Lahti, and J.-P.\ Pellonpää, {\em Covariant localizations in the torus and the phase observables,}
J.\ Math.\ Phys.,  {\bf 43} (2002), pp.\ 693-704.

\bibitem{ChVe} J.\ P.\ R.\ Christensen and J.\ Vesterstr\o m, {\em A note on extreme positive definite matrices,} Math.\ Ann., {\bf 244} (1979), pp.\ 65-68.

\bibitem{DA}G.\ M.\ D'Ariano, {\em Extremal covariant quantum operations and positive operator valued measures,} J.\ Math.\ Phys., {\bf 45} (2004), pp.\ 3620-3635.

\bibitem{GrPiWa} R.\ Grone, S.\ Pierce, and W.\ Watkins, {\em Extremal correlation matrices,} Linear Algebra Appl., {\bf 134} (1990), pp.\ 63-70.

\bibitem{Ho}A.\ S.\ Holevo, {\em Covariant measurements and uncertainty relations,} Rep.\ Math.\ Phys., {\bf 16} (1979), pp.\ 385-400. 

\bibitem{NCQM} J.\ Kiukas, P.\ Lahti, and K.\ Ylinen, {\em Normal covariant quantization maps,} J.\ Math.\ Anal.\ Appl., {\bf 319} (2006), pp.\ 783-801.

\bibitem{Li}B.\ Li, {\em Real operators algebras}, World Scientific, Singapore 2003.

\bibitem{LiTa} C.-K.\ Li and B.-S.\ Tam, {\em A note on extreme correlation matrices,} SIAM J.\ Matrix Anal.\ Appl., {\bf 15} (1994), pp.\ 903-908.

\bibitem{Lo}R.\ Loewy, {\em Extreme points of a convex subset of the cone of positive semidefinite matrices,} Math.\ Ann., {\bf 253} (1980), pp.\ 227-232.

\bibitem{Pa}V.\ Paulsen, {\em Completely bounded maps and operator algebras}, Cambridge Studies in Advanced Mathematics 78, Cambridge University Press, Cambridge 2002. 

\end{thebibliography}
\end{document}